\documentclass{amsart}
\usepackage{latexsym,amsxtra,amscd,ifthen}
\usepackage{amsfonts}
\usepackage{verbatim}
\usepackage{amsmath}
\usepackage{amsthm}
\usepackage{amssymb}

\numberwithin{equation}{section}

\theoremstyle{plain}
\newtheorem{theorem}{Theorem}[section]
\newtheorem{theorema}[theorem]{Theorem}

\newtheorem{prop}[theorem]{Proposition}
\newtheorem{lemma}[theorem]{Lemma}
\newtheorem{cor}[theorem]{Corollary}
\newtheorem{corollarium}[theorem]{Corollary}

\newtheorem{quest}[theorem]{Question}

\theoremstyle{definition}
\newtheorem{defi}[theorem]{Definition}

\newtheorem{definitia}[theorem]{Definition}

\newcommand{\PP}{\ensuremath{\mathbb P}}

\def\la{\langle}
\def\ra{\rangle}
\def\t{\widetilde}

\def\Ann{\mbox{Ann\,}}
\def\ann{\mbox{Ann\,}}

\def\char{\mbox{char\,}}

\def\const{\mbox{Const}}
\def\deg{\mbox{deg\,}}
\def\dim{\mbox{dim\,}}
\def\ker{\mbox{ker\,}}

\def\rate{\mbox{rate\,}}

\def\tor{\mbox{Tor\,}}
\newcommand{\td}{\mathop{\mathrm{tr.~deg.}}\nolimits}

\begin{document}

\title{Linear equations \\
         over noncommutative graded rings%
%\thanks{Partially
%supported by the grant 02-01-00468 of the Russian Basic Research Foundation.}%
}

\author{Dmitri Piontkovski}

      \address{ Central Institute of Economics and Mathematics\\
                      Nakhimovsky prosp. 47, Moscow 117418,  Russia}
\thanks{Partially
supported by the grant 02-01-00468 of the Russian Basic Research Foundation}

\email{piont@mccme.ru}

\begin{abstract}

We call a graded connected algebra $R$ {\it effectively coherent},
if
 for every linear equation over $R$  with homogeneous coefficients of degrees at most
$d$,
 the degrees of generators of its module of solutions are bounded by some function
$D(d)$.
 For commutative polynomial rings, this property has been established by Hermann in 1926.
 We establish the same property for several classes of noncommutative algebras, including
  the most common class of rings in noncommutative projective geometry,
 that is, strongly Noetherian rings, which includes Noetherian PI algebras and Sklyanin
algebras.
 We extensively study so--called {\it universally coherent} algebras, that is,
 such that the function $D(d)$ is bounded by $2d$ for $d \gg 0$.
 For example, finitely presented monomial algebras belong to this class, as well as
 many algebras with finite Groebner basis of relations.
\end{abstract}



\maketitle

\section{Introduction}

\subsection{Overview}

Let $R$ be a connected graded associative algebra over a field $k$.
We discuss the solutions of a linear equation
\begin{equation}
\label{eq}
a_1 x_1 + \dots + a_n x_n = 0
\end{equation}
over $R$, where $a_1,  \dots ,a_n $ are homogeneous elements of $R$
or of a (free) $R$--module  $M$,
and $x_1,  \dots ,x_n$   are indeterminates.
When and how can such an equation be solved, and how does one describe the solutions?

For general finitely generated $R$, there is no algorithm even to check if a solution
exists
(at least in the nonhomogeneous case for a free module $M$, see~\cite{um}).
Also, the set of solutions $\Omega$ may be infinitely generated as a submodule of the
free
module $R^{n}$ (for example, over the algebra $F \otimes F$, where $F$ is a free
associative algebra with a large number of generators).
That is why it seems reasonable to restrict the class of algebras under consideration
to algebras $R$ such that the module of solutions of an equation~(\ref{eq}) is finitely
generated.
Such algebras are called (right)  {\it coherent}~\cite{bur, faith}; this class includes
all Noetherian algebras,
free associative algebras, and many other examples.
However, the condition of coherence does not
in general give a way to find all the generators of $\Omega$:
we can find the generators one by one, but when we have to stop?

If $R$ is a commutative affine algebra, there is an easy (but not the most effective) way
to find the generators of $\Omega$.
It was established by Hermann~\cite{her} in 1926 that there exists a function   $D_R :
{\bf N} \to {\bf N}$
such that $\Omega$ is generated in degrees at most $D_R(d)$  provided that all
coefficients $a_i$ have degrees at most $d$.
So, to find all the solutions, it is sufficient to find the solutions
in the finite--dimensional vector space $R_{\le D_R(d)}$:
it is a standard exercise in linear algebra.

The main object of this paper is to consider non-commutative algebras $R$ which admit
such a function $D_R(d)$.
We call these algebras {\it effectively coherent}.
An analogous concept for commutative (non)graded algebras has been introduced by
Soublin~\cite{ucoh}.
A commutative ring $R$ is called     uniformly   coherent if there is a function
$\Delta_R: {\bf N} \to {\bf N}$
such that $\Omega$ is generated by at most $\Delta_R (n)$ elements. It was shown
in~\cite{coh2}
that an affine or local Noetherian commutative ring is uniformly   coherent if and only
if its dimension is at most two.

Fortunately, our graded analogue of this concept is more
common. We show that most rings considered in
non-commutative projective geometry, that is, strongly
Noetherian algebras~\cite{asz} over algebraically closed
fields, are effectively coherent ($R$ is called (right)
{\it strongly Noetherian} if $R \otimes C$ is right
Noetherian for every commutative Noetherian $k$--algebra
$C$). This class includes, in particular, Sklyanin
algebras, Noetherian PI algebras (in particular, standard
Noetherian semigroup algebras of polynomial
growth~\cite[Theorem~3.1]{gijo}), Noetherian domains of
Gelfand--Kirillov dimension two, and Artin--Shelter regular
algebras of dimension three~\cite{asz}. Also, free
associative algebras and finitely presented monomial
algebras are effectively coherent as well. Every coherent
algebra over a finite field $k$ is effectively coherent;
however, if $k$ contains two algebraically independent
elements or has zero characteristic, then there exist
Noetherian but not effectively coherent algebras.

The class of effectively coherent algebras is closed under
extensions by finitely presented modules, free products,
direct sums, and taking Veronese subalgebras (in the
degree--one--generated case). Every finitely presented
graded module $M$ over such an algebra is effectively
coherent as well, that is, there is a similar function
$D_M(d)$ which bounds the degrees of generators in the case
of the equation~(\ref{eq}) over $M$. This means that
homogeneous linear equations over such modules are
effectively solvable as well.
%Moreover, every system of homogeneous linear equations over a finitely %presented module
%$M$ or over the ring $R$ itself
%is effectively soluble as well as a single equation.

The degree bound function $D_R(d)$ for the commutative
polynomial ring $R$ grows as a double exponent. That is why
we cannot hope that there is a wide class of noncommutative
algebras $R$ with slow growth of  $D_R(d)$. However, there
are interesting classes of algebras with linear growth of
$D_R(d)$. We call an algebra {\it universally coherent} if
$D(d) \le 2d$  for all $d \gg 0$. We investigate this class
of algebras and more general classes, so--called algebras
with Koszul filtrations and with coherent families of
ideals. In particular, it is shown that every finitely
presented module over a universally coherent algebra has
rational Hilbert series, including the algebra $R$ itself,
and $R$ has finite Backelin's rate (that is, there is a
number $r$ such that every space $\tor_i^R(k,k)$ is
concentrated in degrees at most $r i$). Free associative
algebras are  universally coherent, as well as finitely
presented monomial algebras and {\it algebras with
$r$--processing}~\cite{pi2}, that is, algebras with finite
Groebner basis such that the normal form of a product of
two their elements can be calculated by the product of
their normal forms via a bounded number of reductions. The
main property of such algebras is that every right--sided
ideal $I$ has finite Groebner basis: it
 consists of elements
of degree less than $d+r$~\cite[Theorem~5]{pi2}, where $I$ is generated in degrees at most $d$.

Note that the most effective modern method to solve an
equation of type~(\ref{eq}) is based on the theory of
Groebner bases~\cite{iou, pi1, green, pi2}. Let $I$ be a
submodule of $M$ generated by the coefficients $a_1, \dots,
a_n$. In this method, we can calculate the Groebner bases
of relations of $R$, of relations of $M$, and of the
submodule $I$ up  to degree $D(d)$, and then find a
generating set of the relations of degrees at most $D(d)$
between these elements of Groebner basis of $I$, again
using standard Groebner theory methods. The calculation of
the Groebner bases above may be done in the same way as the
usual calculation of Groebner bases of ideals in algebras,
since $I$ is an ideal in the trivial extension algebra $R'
= M \oplus R$; this calculation is equivalent also to
finding the two-sided Groebner basis of relations of the
larger trivial extension $R'/I \oplus R'$ (a similar trick
has been described in~\cite{hey}).

\subsection{Motivation}

Despite of the famous recent progress in noncommutative
projective geometry, no general noncommutative version of
computational methods of algebraic geometry is known. In
this paper, we try to show that a ``computational
noncommutative geometry'' is possible. At least, if $R$ is
a ``ring of noncommutative projective
geometry''~\cite{kel}, then there exist algorithms to solve
linear equations over $R$ (since $R$ is usually strongly
Noetherian), therefore, to calculate the relations and the
minimal projective resolution of a finite module (because
$R$ has often finite global dimension).

\subsection{Notation and assumptions}

We will deal with ${\bf Z}_+$--graded connected associative algebras over a fixed field
$k$, that
is, algebras of the form   $R = \bigoplus_{i \ge 0} R_i$  with $R_0 = k$.
All modules and ideals are graded and right-sided.

A solution $x = (x_1, \dots, x_n )$ of the equation~(\ref{eq}) is called {\it
homogeneous}, if there is $D>0$
(the {\it degree} of $x$) such that $\deg a_i + \deg x_i = D$ for all nonzero $x_i$.

If a sequence  $a = \{ a_1, \dots, a_n \}$ of homogeneous
elements in an $R$--module $M$ generates a submodule $I$,
let $b = \{ a_1, \dots, a_m\}$ ($m \le n$) be a minimal
subsequence of $a$ generating the same submodule. Let
$\Omega^x$  be the module of solutions of the
equation~(\ref{eq}), let $\Omega^y$ be the module of
solutions of the corresponding equation for $b$
\begin{equation}
\label{eq-b}
     a_1 y_1 + \dots + a_m y_m = 0,
\end{equation}
and let $D_x$ and $D_y$ be the maximal degrees of
homogeneous generators of $\Omega^x$ and $\Omega^y$. It is
easy to see that $D_x \le \max \{ D_y ,d \}$. Therefore, we
may (and will) always assume that the coefficients in the
equation~(\ref{eq}) {\it minimally} generate some submodule
$I = a_1 R + \dots + a_n R \subset M$.

  For an $R$--module $M$, we will denote by $H_i M$
the graded vector space $\tor_i^R (M,k)$. By $H_i R$ we
will denote the graded vector space $\tor_i^R (k,k)  = H_i
k_A$. In particular, the vector space $H_1 R$ is isomorphic
to the linear $k$--span of a minimal set of homogeneous
generators of $R$, and  $H_2 R$ is isomorphic to the
$k$--span of a minimal set of its homogeneous relations.
Analogously, the space $H_0 M$ is the span of generators of
$M$, and $H_1 M$ is the span of its relations.

Let $m(M) = m_0(M)$ denote the supremum of degrees of
minimal homogeneous generators of $M$: if $M$ is just a
vector space with the trivial module structure, it is
simply the supremum of degrees of elements of $M$. For
$i\ge 0$, let us also put $m_i (M) := m (H_i M) = \sup \{ j
| \tor_i^R(M,k)_j \ne 0 \}$. Similarly, let us put $m_i(R)
= m (H_i R) = m_i (k_R)$. For example, $m(R) = m_0(R)$ is
the supremum of degrees of the generators of $R$, and $m_1
(R)$ (respectively, $m_1 (M)$) is the supremum of degrees
of the relations of $R$ (resp., of~$M$). In other words, if
a module $I$ is minimally generated by the coefficients
$a_1, \dots, a_n $ of the equation~(\ref{eq}) and $\Omega$
is the module of solutions of this equation, then
$m(\Omega) = m_1(I)$.

Note that the symbols $H_i R$ and  $m_i R$ for an algebra
have different meaning that  the respective symbols $H_i
R_R$ and $m_i R_R$ for $R$ considered as a module over
itself; however, the homologies $H_i R_R$ are trivial, so
that there is no place for confusion.

\begin{defi}
For a finitely generated module $M$,
let us define a function $D_M : {\bf N} \to {\bf N} \cup \{ \infty \} $
by taking $D_M (d) = \sup \{  m_1 (L) | L \subset M, m_0 (L) \le d \}$.
\end{defi}
This means that every submodule  $L \subset M$
 generated in degrees at most $d$
has relations in degrees at most $D = D_M (d)$, and that
the module of solutions of every linear equation~(\ref{eq})
with coefficients of degrees at most $d$ in $M$ is
generated in degrees at most $\max \{ D(d), d \}$.

%For every $R$--module $M$ we define the function $D_M : $

For a graded locally finite vector space (algebra,
module...) $V$, its Hilbert series is defined as the formal
power series $V(z) = \sum_{i \in {\bf Z}}  (\dim V_i) z^i$.
For example, the Euler characteristics of a minimal free
resolution of the trivial module $k_R$ leads to the formula
\begin{equation}
\label{euler}
    R(z)^{-1} = \sum_{i \ge 0} (-1)^i H_i  R(z). \end{equation}
As usual, we write $\sum _{i \ge 0} a_i z^i = o(z^n)$ iff $a_i =0$ for $i \le n$.

Let us introduce a lexicographical total order on the set of all power series with
integer coefficients, i.e., we put
$\sum_{i \ge 0} a_i z^i >_{lex} \sum_{i \ge 0} b_i z^i$ iff
there is $q \ge 0$  such that $a_i = b_i$  for $i < q$ and $a_q > b_q$.
This order extends the coefficient-wise partial order  given by
$\sum_{i \ge 0} a_i z^i \ge \sum_{i \ge 0} b_i z^i$ iff $a_i \ge b_i$ for all $i \ge 0$.

\subsection{Results}

Our technique is based on the investigation of Hilbert
series of algebras and modules. We begin with recalling a
classical theorem of Anick on the Hilbert series of
finitely presented algebras: the set of Hilbert series of
all $n$--generated algebras $R$ with $m_1(R), m_2(R) <
\const$ satisfies the ascending chain condition with
respect to the order~$ >_{lex} $. Then we improve this
theorem for the algebras with additional condition $m_3(R)
< \const$: that is, we state
\begin{theorema}[Theorem~\ref{main}]
\label{main-intro} Given four integers $n,a,b,c$, let
$D(n,a,b,c)$ denote the set of all connected algebras $A$
over a fixed field $k$ with at most $n$ generators such
that $m_1 (A) \le a, m_2(A) \le b$, and $m_3(A) \le c$.
Then the set of Hilbert series of algebras from
$D(n,a,b,c)$ is finite.
\end{theorema}
% the set of Hilbert series of such algebras is finite (Theorem~\ref{main}).
 This
 additional restriction $m_3(R)
< \const$
 (the weakest among all considered in this paper)
is discussed in subsection~\ref{ss m_3}.
 We give also a version of both these theorems for modules
(subsection~\ref{ss mod}):
 in particular, if $m_i(R) < \const$ for $i = 1,2,3$, then the set of Hilbert series of
ideals
 $I \subset R$ with $m_0(I), m_1(I) < \const$ is finite.

 In section~\ref{s EC}, we introduce and study effectively coherent rings.
 First, we give several criteria for a ring to be effectively coherent
 and show that finitely presented extensions,
 free products, free sums, and Veronese subrings of  effectively coherent rings are
 effectively coherent as well, and give appropriate estimates for the function $D(d)$.
 %In the next subsection~\ref{ss SN},
 Futher,
 we introduce other effectivity properties of
graded algebras, related to
 Hilbert series of their finitely presented modules. Let $M$ be a  finitely presented
module, and
 let $L$ run through
 the set of all
 its finitely generated submodules. We say that $M$ is
 {\it effective for generators} (respectively, {\it effective for series}), if,
 given the Hilbert series $L(z)$ (resp., given $m(L)$), there are only finite
 number of possibilities for $m(L)$   (resp., for $L(z)$). An algebra $R$ is said to be
 effective for generators  (resp., for series), if every
finitely presented $R$--module
satisfies this property.
 The relations of these properties to effective coherence are the following:
 if a coherent algebra $R$ is effective both for series and for generators, then it is
effectively coherent,  and every effectively coherent algebra is effective for series.
 Also, we show that the properties of Hilbert series of finitely generated modules over strongly
Noetherian algebras
 established in~\cite[Section~E4]{asz} imply both effectivity for generators and for series:
 in particular, we obtain
\begin{theorema}[Corollary~\ref{str_noet_effect}]
\label{str_noet_effect-intro} Every strongly Noetherian
algebra over an algebraically closed field  is effectively
coherent.
\end{theorema}
 However, there are Noetherian algebras which do not satisfy any of our
 effectivity
properties:
 for example, one of Noetherian but not strongly Noetherian algebras from~\cite{rog}
 (namely, the graded algebra $R_{p,q}$ generated by two Eulerian derivatives introduced
in~\cite{jor}).

 In section~\ref{rate}, we study algebras, not necessary coherent, but having a
lot of finitely presented ideals.
 Such special families of ideals were first introduced for quadratic commutative algebras as
{\it Koszul filtrations}~\cite{crv, ctv};
 then this notion has been generalized to non-quadratic commutative~\cite{cnr}
 and to quadratic noncommutative~\cite{kfx} algebras.
Here we consider the most general version, which is called
{\it coherent family}
of ideals. A family ${\bf F}$ of finitely generated ideals in $R$ is said to be coherent
if $0 \in {\bf F}, R_{\ge 1} \in {\bf F}$, and
 for every  $0 \ne I \in {\bf F}$
    there are $J \in {\bf F}$ and  $ x \in I$
    such that $I \ne J, I = J + x R, m(J) \le m(I)$, and  the ideal
    $(x:J) :=  \{ a \in R | xa \in J \}$ also belongs to ${\bf F}$.
A {\it degree} of ${\bf F}$ is the supremum of degrees of
generators of ideals $I \in {\bf F}$. Coherent families of
degree one are called Koszul filtrations; they do exist in
many commutative quadratic rings (such as coordinate rings
of some common varieties), in algebras with generic
relations, and in quadratic monomial algebras. We show that
if an algebra admits a coherent family ${\bf F}$ of finite
degree, then it has finite Backelin's rate (generalizing an
analogous result in the commutative
case~\cite[Proposition~1.2]{cnr}), and its Hilbert series
is a rational function
% and rational Hilbert series
(generalizing similar result for the algebras with Koszul
filtrations~\cite[Theorem~3.3]{kfx}), and the same is true
for every ideal  $I \in {\bf F}$.

Every ideal in a Koszul filtration is a Koszul module, and
an algebra is coherent if and only if all its finitely
generated ideals form a coherent family. If all ideals of
an algebra $R$, generated in degrees at most $d$,  form a
coherent family, we call $R$ universally $d$--coherent. An
algebra $R$ is called {\it universally coherent}  if it is
universally $d$--coherent for all $d \gg 0$. In fact, these
properties are the properties of the function $D_R(d)$:
\begin{theorema}[Proposition~\ref{D<t+d}, Corollary~\ref{D<2d}]
\label{kf-intro} Let $R$ be a finitely generated graded
algebra.

(a) $R$ is universally $d$--coherent iff $D_R(t) \le t+d$
for all $t \le d$.

(b) $R$ is universally coherent iff $D_R(d) \le 2d$ for all
$d \gg 0$.

In particular, any universally coherent algebra is
effectively coherent.
\end{theorema}
Commutative 1--universally coherent algebras are called
{\it universally Koszul}; they has been studied
in~\cite{conca, con2}. We show that the $d$--th Veronese
subring of a generated in degree one universally
$d$--coherent algebra $R$ is universally Koszul, therefore,
such an algebra is (up to a shift of grading) a Koszul
module over a universally Koszul algebra.

Some noncommutative examples of universally coherent
algebras are considered in subsection~\ref{perera}, that
is,
%these are
finitely presented monomial algebras and, more
generally, a class of algebras with a finite Groebner basis
of relations (algebras with $r$-- processing), which were
introduced in~\cite{pi2}.

\subsection{Acknowledgment}

I am grateful to to Alexander Kuznetsov for useful remarks
and to Daniel Rogalski for helpful conversations. Also, I
am grateful for hospitality to Mittag--Leffler Institut
where this work has been accomplished. Finally, I thank the
anonymous referee for wide and detailed comments.

\section{Sets of Hilbert series}

\subsection{The condition $\dim \tor_3^R(k,k) < \infty$}

\label{ss m_3}

In this section, we sometimes consider finitely presented
algebras $R$
such that
 %with extra condition
 $m_3(R)  <
\infty$. Before studying their Hilbert series, let us say a
few words about this inequality.

First, all common finitely presented algebras (such as
Noetherian, coherent, Koszul etc) do satisfy this
condition. In fact, it is the weakest restriction on $R$
among all that are considering in this paper. In a coherent
ring, the module of solutions of any linear equation over a
free module is finitely generated; in general, there is a
particular linear equation in a free module which has
finite basis of solutions if and only if $m_3(R)  <
\infty$.

Indeed, let $a = \{ a_1, \dots, a_g\}$ be a minimal set of
homogeneous generators of $R$, and let $f = \{ f_1, \dots,
f_r\}$ be a minimal set of its   homogeneous relations. Let
$f_j = \sum_{i=1}^g a_i b_j^i$ for $j = 1 \dots r$. In the
minimal free resolution of $k_R$
$$
     \dot \to H_3(R) \otimes R   \to H_2(R) \otimes R   \to H_1(R) \otimes R  \to R \to k
\to 0
$$
we see that $H_1(R)$ is the span of $a$, and $H_2(R) $ is
the span of $f$. Let $\t f_j = \sum_{i=1}^g a_i \otimes
b_j^i \in ka \otimes R$ be the image of $f_j \otimes 1$ in
the free module $M  = ka \otimes R$. Consider the following
equation with coefficients in $M$:
$$
      \t f_1 x_1  + \dots + \t f_r x_r = 0.
$$
Since the resolution above is minimal, every minimal space of generators of the solution
module $\Omega$  of this equation is isomorphic to $H_3(R)$.

However, in general, given a presentation $(a, f)$ of an
algebra $R$,  there does not exist an algorithm to decide
if the condition $m_3(R) < \infty$ holds. This has been
shown in~\cite{an-undec} for Roos algebras, that is,
universal enveloping algebras of quadratic graded Lie
superalgebras.

\subsection{Hilbert series of finitely presented algebras}

The following well--known theorem describes an interesting
property of Hilbert series of finitely presented algebras.

\begin{theorema}[{\cite[Theorem~4.3]{an4}}]
\label{ant} Given three integers $n,a,b$, let $C(n,a,b)$ be
the set of all $n$--generated connected algebras $R$ with
$m_1(R) \le a$ and $m_2(R) \le b$ and let ${\mathcal H}
(n,a,b)$ be the set of Hilbert series of such algebras.
Then the ordered set $({\mathcal H} (n,a,b) ,  >_{lex} )$
admits no infinite ascending chains.
\end{theorema}

The example of an infinite {\it descending} chain of
Hilbert series in the set $C(7,1,2)$ is constructed in
~\cite[Example~7.7]{an4}. All algebras in this chain have
global dimension three, but in the vector spaces $H_3 R$
there are elements of arbitrary high degree.

The following theorem shows, in particular, that the last property is essential for such
examples.

\begin{theorema}
\label{main} Given four integers $n,a,b,c$, let
$D(n,a,b,c)$ denote the set of all connected algebras $A$
over a fixed field $k$ with at most $n$ generators such
that $m_1 (A) \le a, m_2(A) \le b$, and $m_3(A) \le c$.
%of degree at most $a$, whose relations are
%concentrated in degrees at most $b$ such that the
%graded vector space $\tor_3^A(k,k)$  is  concentrated in degrees at most $c$.
Then the set of Hilbert series of algebras from
$D(n,a,b,c)$ is finite.
\end{theorema}

For the set $D(n,1,2,3)$,  Theorem~\ref{main} was proved
in~\cite[Corollary~2 in subsection~4.2 and Remark~1 after
it]{pp} (in a different way, using a geometrical
technique). In particular, given a number $n$, the set of
Hilbert series of $n$--generated quadratic Koszul algebras
is finite. For the class of degree--one generated algebras
of bounded Backelin's rate (that is, when there is a number
$r$ such that for every $i$ the vector space $H_i R$ is
concentrated in degrees at
most $r i$), a similar statement
has been proved by L.~Positselski (unpublished). %see~\cite{kfx}).

We need the following standard version of Koenig lemma.
\begin{lemma}
\label{Koenig}
Let $P$ be a totally ordered set satisfying both ACC and DCC. Then $P$ is finite.
\end{lemma}

%\begin{proof}[Proof of Lemma]
%The condition ACC (DCC) means that every subset of
%$P$ contains its supremum (infinum).
%Let $P_M$  be a subset of all elements $x$ of $P$ such that the set
%$x_> = \{ y | y  >x\}$ is finite.
%By ACC, $P_M$ is non-empty. By DCC, there is $\inf P_M \in P_M$,
% hence $P_M$ is finite.
%Suppose that
%the set $Q = P \setminus P_M $ is non-empty. Then there is $s = \sup Q \in Q$.
%Therefore, the set
%$s_> = P_M$ is finite, in contradiction to the assumption $s \in Q$. %\end{proof}

\begin{proof}[Proof of Theorem~\ref{main}]
Consider a connected algebra $A$ with a minimal space of
generators $V$  and a minimal space of relations $R \subset
T(V)$. Choose a homogeneous basis $f = \{ f_i\}$ in $R$.
Let $I$ be the ideal in $ T(V)$ generated by $R$, and let
$G$ be the graded algebra associated to the $I$--adic
filtration on $ T(V)$. By~\cite[Theorem~3.2]{gd3}, we have
$H_j  G = H_j A \oplus H_{j+1} A$ for all $j \ge 1$; in
particular, the space of generators of $G$ is isomorphic to
$V \oplus R$. Let $\tilde f = \{ \tilde f_i\}$ be the set
of generators of the second summand, corresponding to the
basis $f = \{ f_i\}$ in $R$.

Note that $G$ is generated by its subsets $A$ and $\tilde f
$, and its grading extends the grading of $A$ if $\deg
\tilde f_i = \deg f_i$.  Moreover, the algebra $G$ is the
quotient of the free product $A * k \la \tilde f \ra $ by
the ideal generated by some elements of $A \otimes \tilde f
\otimes A$~\cite[Section~3]{gd3}, so that we can consider
another of $G$ given by $\deg' a = \deg a$ for $a \in A$
and $\deg' \tilde f_i = \deg f_i - 1$.
%By the above
%properties of relations of $G$, it remains graded.
Let us denote the same algebra $G$ with this new grading by
$C = C(A)$. It follows from the consideration
in~\cite[proof of Lemma~5.5]{gd3} that its homology groups
 are given  by the formula $H_j  C = H_j A \oplus H_{j+1}
A [1]$ for all $j \ge 1$. By the formula~(\ref{euler}), we
have
$$
      C (z)^{-1} = \sum_{i \ge 0} (-1)^i H_i C (z) =
$$
$$
     = \sum_{i \ge 0} (-1)^i H_i A (z) + z^{-1}
\left( 1-V(z) -
     \sum_{i \ge 0} (-1)^i H_i A (z)     \right) =
$$
$$
    = A (z)^{-1} (1 - z^{-1}) + z^{-1} (1-V(z)).
$$

Now suppose that two connected algebras $A,B$ have the same
graded vector space of generators $V$. Suppose $A(z)
>_{lex}  B(z)$, that is, $A(z) - B(z) = p z^q + o (z^q)$
for some $p> 0, q> 1$. We have
$$
C (A) (z) - C (B)(z) = (C (B)(z)^{-1} - C (A) (z)^{-1}) C (A) (z) C (B) (z) =
$$
$$
= (B(z)^{-1} - A (z)^{-1}) (1 - z^{-1}) (1 + o(1)) =
    \dfrac{A (z) - B(z)}{z A(z) B(z)} (1 + o(1)) =
$$
$$
 = - p z^{q-1}  + o (z^{q-1}),
$$
that is, $C (A) (z) <_{lex} C (B)(z)$.

Now, we are ready to prove the theorem.
Assume (ad absurdum) that the set  $D(n,a,b,c)$ is infinite for some $n,a,b,c$. By
Lemma~\ref{Koenig},
there is an infinite descending chain     $Ch_0$
$$
A^{(1)} (z) >_{lex} A^{(2)} (z) >_{lex} \dots
$$
Since for algebras in the set $D(n,a,b,c)$ there are only
finite number of possibilities for the number and degrees
of generators, the chain   $Ch_0$    contains an infinite
subchain
  $Ch_1$
$$
A^{1} (z) >_{lex} A^{2} (z) >_{lex} \dots,
$$
where all  algebras $A^i$ are generated by the same graded
vector space $V$. It follows that we have an accending
chain $C (Ch_1)$:
$$
     C (A^{1} )  (z) <_{lex}   C (A^{2} )  (z) <_{lex} \dots
$$
For every algebra $C(A^i)$, its generators are concentrated
in degrees at most $a' = \max (a, b-1)$ and  relations  are
concentrated in degrees at most $b' = \max (b, c-1)$.
Moreover, for the number of its generators we have the
following estimate:
$$
    \dim H_1 C(A^i) = \dim H_1 A^i +  \dim H_2 A^i \le a + a^b =: n'.
$$
We deduce that an infinite ascending chain $C (Ch_1)$
consists of algebras from the set $C(n',a',b')$, in
contradiction to Theorem~\ref{ant}.
\end{proof}

\subsection{Modules and ideals}

\label{ss mod}

The following gives module versions of Theorems~\ref{ant}
and~\ref{main} for modules.

\begin{prop}
\label{mod-gen}

Let $n,a,b,c, m, p_1, p_2, q, r$  be 9 integers.

(a) Let $R$ be an algebra from $C(n,a,b)$, and let $CM =
CM(m, p_1, p_2, q)$ denote the set of all graded right
$R$--modules $M$ with at most $m$ generators such that $M_i
= 0$ for $i<p_1$ and $i > p_2$ (in particular, $m_0(M) \le
p_2$), and $m_1(M) \le q$. Then the ordered set of Hilbert
series of modules from $CM$ satisfies ACC.

(b)   Let $DM =DM(n,a,b,c, m, p_1, p_2, q, r)$ denote the
set of all graded right modules  over algebras from
$D(n,a,b,c)$ with at most $m$ generators such that  $M_i =
0$ for $i<p_1$, $m_0(M) \le p_2, m_1(M) \le q$, and $m_2(M)
\le r$. Then the set of Hilbert series of modules from $DM$
is finite.
\end{prop}

\begin{proof}
If $p_1 < 0$, let us shift the grading of all modules by
$1-p_1$ and consider the sets of Hilbert series of the
modules from $ CM(m, 1, p_2-p_1+1, q-p_1+1)$ (respectively,
$DM(n,a,b,c, m, 1, p_2-p_1+1, q-p_1+1, r-p_1+1)$). Since
these new sets are in bijections with $CM$ and $DM$, we may
assume that $p_1
>0$, that is, that all our modules are generated in
strictly positive degrees.

Let $R$ be an algebra in $C(n,a,b)$ (respectively, in
$D(n,a,b,c)$), and let $M$ be an $R$--module contained in
$CM$ (respectively, in $DM$). Consider its trivial
extension $C_M = M \oplus R$. By the classical
formula~\cite{trivext0} for Poincare series of trivial
extensions, we have
$$
    P_{C_M} (s,t) = \frac{P_R(s,t)}{1-s P_M (s,t)}
        = P_R(s,t) (1+s P_M (s,t) + s^2 P_M (s,t)^2 +\dots)
$$
(where $P_{-}(s,t) = \sum_{i \ge 0} s^i H_i(-)(t)$),
hence
$C_M \in CM(N,A,B)$  (resp., $C_M \in D(N,A,B,C)$)
 for some $N,A,B$(,$C$) depending on $n,a,b, m, p_1, p_2, q $($,c,r$).
 In the case $(a)$, we apply Theorem~\ref{ant} and conclude that the set of Hilbert
series of such algebras
 $C_M$ satisfies ACC; in the case $(b)$, we also apply Theorem~\ref{main} and find that
there is only
 a finite number of possibilities for $C_M (z)$.
Since $M(z) = C_M (z) - R(z)$, the same is true for the set of Hilbert series $M(z)$.
\end{proof}

\begin{cor}
\label{mod}
Let $D > 0$ be an integer.

(a)
If $R$ is a connected finitely presented algebra,
then the set of Hilbert series of right-sided ideals
in $R$ generated in degree at most $D$ satisfies DCC.

(b) If, in addition, $m_3(R) < \infty $, then
the set of Hilbert series of right-sided ideals
in $R$ having generators and relations in degrees at most $D$
is finite.
\end{cor}

\begin{proof}
If $I$ is an ideal in $R$ and $M =R/I$, then the exact sequence
$$              0 \to I \to R \to M \to 0
$$
implies that $I(z) = R(z) - M(z)$ and $m_{i+1}(M) =
m_i(I)$. It remains to apply Proposition~\ref{mod-gen} to
the set of such modules $M$.
\end{proof}

There is a class of algebras for which the finiteness of
the set of Hilbert series of right ideals can be proved
without the assumption on degrees of relations.
% For some algebras $R$, the additional assumption that the
%relations of the ideals must have degrees at most $D$ may
%be omitted in the case $(b)$.
Such algebras will be
considered in the next section.

\section{Effective coherence}

\label{s EC}

\subsection{Effectively coherent rings}

All algebras below are connected graded, all modules (and ideals) are right and graded,
as before.

Recall that an algebra is called (graded or projective)
coherent if every its finitely generated ideal is finitely
presented, or, equivalently, the kernel of any map  $M \to
N$ of two finitely generated free modules is finitely
generated. Other equivalent conditions may be found
in~\cite{faith,bur,chase}. In particular, every Noetherian
ring is coherent, while the free associative algebras are
coherent but not Noetherian.

Let us give an effective version of this definition.

\begin{defi}
\label{def-coher} Let $A$ be an algebra. A finitely
generated $A$--module $M$ is called  {\it effectively
coherent} if $D_M(d)$ is finite for every $d > 0$.

The algebra $A$ is called   {\it effectively coherent}, if
it satisfies either of the following equivalent conditions:

(i) $A$ is  effectively coherent as a module over itself;

(ii) every finitely presented $A$--module is  effectively coherent;

(iii) for every finitely presented $A$--module $M$ there is a sequence of functions
$\{ D_i : {\bf N} \to {\bf N} \}$ such that, whenever a submodule $L \subset M$
is generated in degrees at most $d$, the graded vector spaces $Tor_i^A(L, k)$
are concentrated in degrees at most $D_i (d)$ for all $i \ge 0$.
\end{defi}

\begin{proof}[Proof of equivalence]

We begin with
\begin{lemma}
\label{EC-triple}
Let $A$ be an algebra. In an exact  triple of $A$--modules
$$
      0 \to K \to M \to N \to 0,
$$
if any two of these three modules are effectively coherent,
then the third is.

In this case
$$
     D_K(d) \le  D_M(d) \le D_K( \max\{  d, D_N(d) \} )
$$
and
$$
           D_N(d) \le D_M( \max\{  d, m(K) \} ).
$$
\end{lemma}

This  Lemma can be shown in the same way as an analogous
statement for coherent modules, but involving appropriate
estimates for the functions $D_M, D_N, D_K$. Following
N.~Bourbaki~\cite[Exercise~10 to~\S3]{bur}, we leave the
proof to the reader.

Let us return to the proof of equivalence in
Definition~\ref{def-coher}. Since every finite free module
of rank greater than one is a direct sum of free modules of
smaller ranks, every such free module is effectively
coherent provided that $A_A$ is, as well as every its
finite submodule. In particular, a finite presentation
$$
    F'' \stackrel{\phi}{\longrightarrow} F' \to M \to 0
$$
gives the following exact sequence
$$
     0 \to \ker \phi \to   F' \to M \to 0
$$
with effectively coherent first two terms.
This proves $(i) \Longrightarrow (ii)$.

To prove  $(ii) \Longrightarrow (iii)$, just show by induction
that $i$--th syzygy module is effectively coherent.

It remains to point out that the implication $(iii) \Longrightarrow (i)$ is trivial.
\end{proof}

A Noetherian  effectively coherent  algebra is called {\it
effectively Noetherian}: in such algebras, all finite
(=finitely generated) modules are effectively coherent.
Note that the same property of commutative affine algebras
is well-known at least since 1926~\cite{her} when the first
(double--exponential) bound for $D(d)$  for polynomial
rings was found. It is a particular case of effective
Nullstellensatz and effective division problem,
 and there are many papers
(MathSciNet gives about 50)
concerning syzygy degree bounds  $D(d)$
and Betti number degree bounds  $D_i (d)$ for ideals in commutative affine algebras.

We will see in the last subsection that every finitely
presented monomial algebra  is effectively coherent, as
well as coherent algebras with finite Groebner bases
introduced in~\cite{pi2}. A class of  effectively
Noetherian algebras (over an algebraically closed field)
includes so-called strongly Noetherian algebras (in
particular, Noetherian PI--algebras and 3--dimensional
Sklyanin algebras), as we will show in the next subsection.

%We can re--formulate the definition of effective
%coherence in terms of linear equations.
%The next criterion shows that,
%over an  effectively coherent ring, every system %of linear equations is
%effectively solvable as well as a single equation.
%\begin{prop}
%\label{system}
%Let $R$ be an  effectively coherent algebra.
%Then  for every finitely presented
%$A$--module $M$ there is a function
%$\{ D^{sys}_M : {\bf N}^2 \to {\bf N} \}$
%such that, for every system of linear
%equations of type~(\ref{eq})
%of $n$ variables $x_1, \dots, x_n$
% with coefficients in $M$ of degrees at most $d$, its
%module of solutions is generated in degrees at most $D^{sys}_M (n,d)$.
%\end{prop}

Several methods to construct coherent algebras work as well
for effectively coherent ones.

The next criterion is a variation of the classical
criterion of coherence~\cite{bur, chase}.

\begin{prop}
\label{cap-ann}
An algebra $R$ is effectively coherent if and only if
there are two functions $D^\cap , D^{\ann} : {\bf N} \to {\bf N} $
such that for every $a \in R$ we have $m(\ann_R a) \le D^{\ann}(\deg a)$
and for every two right sided ideals $I,J$ with $m(I) \le d, m(J) \le d$ we have
$m(I \cap J) \le D^\cap (d)$.

In this case we have
$$
D^{\ann}(d) \le D(d),   D^\cap (d) \le     \max \{  d, D(d)  \}, \mbox{ and } D(d) \le
\max \{  D^{\ann}(d),    D^\cap (d)\}.
$$
\end{prop}

We begin with
\begin{lemma}
\label{inters}
Let $M$ be a graded module, and
let $K,L$ be two its submodules generated in degrees at most $d$.
Then
$$
         m(K \cap L) \le \max \{ d, D_M(d) \} .
$$
\end{lemma}

\begin{proof}

Consider the following exact triple
$$
    0 \to K\cap L \to K \oplus L \to   K + L \to 0.
$$
The exact sequence of $\tor$'s gives:
$$
    \dots \to \tor_1^R( K + L , k)  \to   \tor_0^R( K\cap L , k)  \to   \tor_0^R(K \oplus
L  , k)  \to \dots
$$
Thus
$$
      m(K \cap L) \le      \max \{  m_1( K + L), m(K), m(L)    \}   \le \max \{  d,
D_M(d)  \} .
$$
\end{proof}

\begin{proof}[Proof of Proposition~\ref{cap-ann}]

The  ``only if'' part follows from Lemma~\ref{inters}. To
show the ``if'' part, we will show that $D(d) \le \max\{
D^\cap  (d), D^{\ann} (d) \}$.
 Let us proceed by induction in
the number of generators  $t$ of an ideal $K \subset A$
with $m(K) \le d $. For $t = 1$, we have $m_1(K) \le
D^{\ann}(d)$. For $t > 1$, we may assume that $ K =   I +
J$, where $I$  and $J$ are generated by at most $(t-1)$
elements. By exact triple
$$
    0 \to I \cap J \to I \oplus J \to   I + J \to 0,
$$
we have
$$
    m_1(K) \le  \max \{  m_1 (I), m_1(J), m (I \cap J) \} .
$$
Here $m (I \cap J) \le D^\cap (d)$ by definition and $m_1
(I), m_1(J) \le \max\{ D^\cap  (d), D^{\ann} (d) \}$ by
induction hypothesis, so the claim follows.
\end{proof}

The following claim is standard for coherent algebras
(see~\cite[Proposition~10]{aberg},\cite[Proposition~1.3]{po}
for two its generalizations).

\begin{prop}
\label{map-coh}
Let $A \to B$ be a map of connected algebras. Suppose that $A$ is (effectively) coherent
and the module $B_A$ is finitely presented. Then $B$ is  (effectively) coherent.
\end{prop}

\begin{proof}

Let $b= m(B_A)$, let $J$ be a right sided ideal in $B$ with
$m(J) = d$, and let $0 \to K \to F \to J \to 0$ be its
minimal presentation with a free $B$--module $F$. Here
$m(F_A) \le m(F_B) +m(B_A) = d+ b$ and
 $m(J_A) \le d +p$, hence $m_1(J_A) \le D_{B_A} (d + p)$.

From the exact sequence of $\tor$'s we have
$$
    \tor_1^A (J,k) \to \tor_0^A (K, k) \to \tor_0^A (F,k).
$$
Therefore,
$$
m_1(J) = m(K_B) \le m(K_A) \le \max \{ m(F_A), m_1(J_A) \} \le  \max \{  d + p ,
D_{B_A} (d + p) \} .
$$
\end{proof}

\begin{cor}
A singular extension of an (effectively) coherent algebra along a finitely presented
module is
(effectively) coherent.
\end{cor}

\begin{prop}
Let $A$ and $B$ be two  (effectively) coherent algebras.
Then their direct sum with common unit $A \oplus B$ and their free product $A * B$
are (effectively) coherent as well.
\end{prop}

The coherence of the free product of two coherent algebras
has been proved in~\cite[Theorem~2.1]{freeprod} (and, in
more general settings, in~\cite[Theorem~12]{aberg}).

\begin{proof}

Let $C = A \oplus B$ and $ E =A * B$. Let $X_A$ and $ X_B$
be minimal homogeneous sets of generators of  $A$ and $B$.

Let $I = \sum_{i=1}^s y_i C$  be a finitely generated ideal
in $C$, where $y_i = a_i + b_i, a_i \in A, b_i \in B$ with
$\deg y_i \le d$. Then $I$ is the factor module of the free
module $F$ with generators $\t y_1, \dots, \t y_s$ by a
syzygy submodule $K$. Then an element $w = \sum_{i=1}^s \t
y_i (\alpha_i + \beta_i) \in F$ belongs to $K$ if and only
if $\sum_{i=1}^s a_i \alpha_i  = 0$  and $\sum_{i=1}^s b_i
\beta_i = 0$. Let $K_A, K_B$  be the syzygy modules of the
ideals in $A$ and $B$ generated by $a_1, \dots , a_s$ and
$b_1, \dots , b_s$, respectively, and let $R_A, R_B$ be
minimal sets of generators of these syzygy modules. The
modules $K_A$ and $K_B$ are submodules of $F$, and $K$ is
equal to the intersection of submodules generated  by $R_A
+ X_B$ and $R_B + X_A$. Hence $K$ is generated by $R_A \cup
R_B$, thus $m_1(I) = \max \{ m(K_A), m(K_B)\} \le \max \{
D_A(d), D_B(d) \} $. So, $C$ is (effectively) coherent if
and only if both $A$ and $B$ are.

Now, let $J$ be a right sided ideal in $E$ with $m(J) = d$, and let $M = E/J$.
The Mayer-Vietoris long sequence~\cite[Theorem~6]{mv}
$$
    \dots \to  \tor_p^k(M,k) \to      \tor_p^A(M,k) \oplus   \tor_p^B(M,k)   \to
\tor_p^E(M,k) \to   \tor_{p-1}^k(M,k) \to \dots
$$
gives an isomorphism
$$
        0 \to   \tor_2^A(M,k) \oplus   \tor_2^B(M,k)  \to   \tor_2^E(M,k) \to 0.
$$
So, $m_1(J) = m_2(M) = \max \{ m_2(M_A), m_2(M_B) \} =  \max \{ m_1(J_A), m_1(J_B) \} $.

For instance, let us estimate $m_1(J_A)$. The set $R$ of
generators of $J$ lies in $E_{\le d}$. Let $J^A$ be the
span of all elements of the form $ux$, where $u \in J, x
\in X_A$, and $\deg u \le d$ but $\deg ux > d$. The subset
$J^B$ is defined in the same way by replacing $A$ by $B$.
Then we have a direct sum decomposition $J = J' \oplus
J''$, where $J' = J_{\le d} \oplus J^A A = J_{\le d} A$
 and $J'' = J^A E X_B E \oplus J^B E$.
Here both summands are $A$--modules, moreover, the second summand $J''$ is a free
$A$--module.
That is why the relations between  the elements of $R \subset J'$
 are the same as the relations of the submodule $J'$.
But $J'$ is a submodule of a finitely generated
free $A$--module $P = E_{\le d}A $, so that
 $     m_1(J_A)  = m_1 (J') \le D_P (d) $.
Involving the analogous estimate for $m_1(J_B)$, we have finally
$$
        m_1 (J) \le  \max \{  D_{E_{\le d}A} (d) ,  D_{E_{\le d}B} (d)   \}  .
$$
\end{proof}

The following corollary will be generalized later in subsection~\ref{perera}.
\begin{cor}
Any free algebra with finitely many generators is effectively coherent.
\end{cor}

The following criterion of coherence of Veronese subrings
has been proved by Polishchuk~\cite[Proposition~2.6]{po}.
Here we give its effective version.

\begin{prop}
\label{DVerona}
Let $A$ be a connected finitely presented algebra generated in degree one.
For every $n \ge 2$, the Veronese subring $A^{(n)} = \bigoplus_{i \ge 0} A_{in}$
is (effectively) coherent if and only if $A$ is (effectively)   coherent.

In the notation of  Proposition~\ref{cap-ann}, we have
$$
          D^{\cap}_{A^{(n)}} (d)  \le  D^{\cap}_A (d) + n-1,     D^{\ann}_{A^{(n)}} (d)
\le    D^{\ann}_A( d) + n-1.
$$
\end{prop}

\begin{proof}

By~\cite[Proposition~2.6]{po}, $A$ is a finitely presented
module over $A^{(n)}$. By Proposition~\ref{map-coh}, it
follows that if $A^{(n)}$  is (effectively) coherent then
so is $A$.

Assume that $A$  is effectively coherent. To show that
$A^{(n)}$ is also effectively coherent, we are going to
apply Proposition~\ref{cap-ann}.

Let  $x \in A^{(n)}$ be an element of degree $d$.
If $\ann_A (x) = \sum_{i=1}^s y_i A$ with $\deg y_i = nq_i - r_i $, where
$0 \le r_i < n$, then  $\ann_{A^{(n)}} (x) = \sum_{i=1}^s y_i A_{r_i} A^{(n)}$,
hence  $m_1(x A^{(n)} ) \le D^{\ann}_A(d) + n-1$.
Therefore, $D^{\ann}_{A^{(n)}} (d) \le D^{\ann}_A(d) + n-1 < \infty$.

Now, let $I= \sum_{i=1}^u a_i A^{(n)}$ and $J= \sum_{i=1}^v b_i A^{(n)}$ be two ideals in
$A^{(n)}$
generated in degrees at most $d$,
and let $K = (\sum_{i=1}^u a_i A)  \cap (\sum_{i=1}^v b_i A) \subset A$.
Let $K = \sum_{i=1}^w c_i A$, where $\deg c_i = nq_i' - r_i' \le D_A(d)$  with $0 \le
r_i' < n$.
Then $I \cap J = K^{(n)} = \sum_{i=1}^w c_i A_{r_i'} A^{(n)}$,
hence $m(I \cap J) \le D_A^\cap (d) + n -1$.
Thus $D^{\cap}_{A^{(n)}} (d) \le D^{\cap}_A (d) + n -1 < \infty$.
\end{proof}

\subsection {Strongly Noetherian algebras are effectively Noetherian}

\label{ss SN}

Effective coherence implies some properties of Hilbert
series. For example, we will see in section~\ref{rate}
below that if $D_A(d) \le 2d$ for $d \gg 0$, then the
Hilbert series $A(z)$  is a rational function. Let us
introduce three other properties of Hilbert series which
are closely related to  effective coherence: the first two
are dual to each other, but the third is stronger.

\begin{defi}
\label{ES-AZ} Let $A$ be an algebra, $M$ be a finite
module, and $L$ runs through the set of its finitely
generated submodules.

{\bf (1)}  $M$ is said to be effective for series if,
given  $m(L)$, there are only
finite number of possibilities for Hilbert series $L(z)$.

{\bf (2)}   $M$ is said to be effective for generators if,
given a Hilbert series $L(z)$, there are only finite number
of possibilities for $m(L)$.

{\bf (3)}  $M$ is said to be Artin--Zhang if for every
formal power series $h(z)$ there is $d > 0$ such that, if
for every $L$ with $L(z) = h(z)$ we have $m(L) \le d$, and
for every $L$ with $m(L) \le d$ and $L(z) = h(z) + o(z^d)$
we  have  $L(z) = h(z)$.

A finitely presented algebra $A$ is said to be effective for series (respectively,
effective for generators, Artin--Zhang),
if every finitely presented $A$--module satisfies such a property.
\end{defi}

For commutative algebras, the effectivity for series
follows from the fact that every ideal (or submodule) $I$
has a Groebner basis of bounded degree. Similar property
has been established also for several classes of ideals in
noncommutative rings (for example, for torsion free finite
modules over 3--dimensional quadratic Artin--Shelter
regular algebras~\cite[Theorem~A]{dnvdb}). The ``effective
for generators'' property is an obvious part of the
Artin--Zhang condition; the property of being effective for
series is naturally dual to  effectivity for generators.

The Artin--Zhang property itself first appeared
in~\cite{az} in the following context. A connected algebra
$A$ is said to be (right) {\it strongly Noetherian} if $A
\otimes R$ is  (right) Noetherian for every Noetherian
commutative $k$--algebra $R$~\cite{asz}. In particular,
Noetherian affine PI algebras, Sklyanin algebras,
Noetherian domains of Gelfand--Kirillov dimension  2,
Artin--Shelter regular algebras of global dimension three,
and some twisted homogeneous coordinate rings
 are strongly Noetherian~\cite{asz}.

\begin{theorem}[{\cite[section~E4]{az}}]
\label{artin-zhang}
Let $A$ be a
strongly Noetherian algebra over an algebraically
closed field $k$. Then $A$ is  Artin--Zhang.
\end{theorem}

A partial converse is given by

\begin{prop}
Every finitely generated Artin--Zhang module is Noetherian.
In particular, every Artin--Zhang algebra is Noetherian.
\end{prop}

The proof consists of the following two lemmas.

\begin{lemma}
\label{ACC_for_AZ} Let $M$ be a finitely generated
Artin--Zhang module. Then the set ${\bf H} = {\bf H}^M$ of
Hilbert series of finitely generated submodules of $M$
satisfies ACC with respect to  the order~``$<_{lex}$''.
\end{lemma}

Notice that for the partial order~``$<$'', the same is
proved in~\cite[Corollary~E4.13]{az}.

\begin{proof}

Assume the converse, that is, that there is an infinite  sequence
$$
    L^1(z ) <_{lex} L^2 (z ) <_{lex}  \dots
$$
of Hilbert series of finite submodules of $M$. Since
$L^i(z) < M(z)$, there exists $\lim_{i \to \infty} L^i(z) =
h(z)$. Let $d$ be the same as in the definition of
Artin--Zhang module. For $i \gg 0$ we have $L^i(z) = h(z) +
o(z^d)$. Let $N^i$ be a submodule of $L^i$ generated by all
elements of $L^i$ having degree at most $d$, i.~e., $N^i =
L^i_{\le d} A$. By Artin--Zhang condition, for $i \gg 0$ we
have $N^i(z) = h(z)$ because $N^i(z) =h(z) + o(z^d)$. Then
we have $h(z) = N^i (z) \le_{lex} L_i(z) <_{lex} h(z)$ for
$i \gg 0$, a contradiction.
\end{proof}

\begin{lemma}
\label{ACC_imply_Noet} Let $M$ be a finitely generated
module such that  $( {\bf H}^M, <_{lex} )$ satisfies ACC.
Then $M$ is Noetherian.
\end{lemma}

\begin{proof}
Assume that $L \subset M$ is minimally generated by an
infinite sequence $x_1, x_2, \dots$, and let $L^i$ denote
its finitely generated submodule $x_1 A + \dots + x_i A$.
Then we obtain an infinite chain in ${\bf H}^M$:
$$
      L^1(z) < L^2(z) < \dots
$$
\end{proof}

\begin{quest}
\label{AZ_ergo_SN?}
Is every Artin--Zhang algebra  over an algebraically
closed field strongly Noetherian?
\end{quest}

Our next purpose is to prove that every strongly Noetherian
algebra over an algebraically closed field is also
effectively Noetherian. To do this, we establish the
following relations between our effectivity properties.

\begin{theorem}
\label{main_effect} Let $A$ be a coherent algebra. Consider
the following properties:

{\sf (AZ)}  $A$ is Artin--Zhang;

{\sf (ES+EG)} $A$ is both effective for series and
effective for generators;

{\sf (EC)} $A$ is effectively coherent;

{\sf (ES)} $A$ is effective for series.

Then there are implications:
$$
    {\sf (AZ)}  \Longrightarrow   {\sf (ES+EG)}  \Longrightarrow   {\sf (EC)}
\Longrightarrow    {\sf (ES)}
$$
\end{theorem}
Note that the assumption that $A$ is coherent is actually
used only in the implication ${\sf (ES+EG)} \Longrightarrow
{\sf (EC)} $.

\begin{proof}
Fix a finitely presented $A$--module $P$. For $d \in {\bf
Z}$, let  ${\bf H}_d$ be the set of Hilbert series of its
submodules generated in degrees at most $d$, and let ${\bf
H} = \bigcup_d {\bf H}_d$.

 ${\sf (AZ)}  \Longrightarrow    {\sf (ES+EG)}$.
Assume that $P$ is Artin--Zhang, and let $h(z)$ be a formal
power series. Let $d = d(h)$ be as in
Definition~\ref{ES-AZ},~${\bf (3)}$.

If $L(z) = h(z)$, then $L_{\le d} A (z)= L(z)$, hence $m(L) \le d$. So, $P$ is  effective
for  generators.

 For any $m \ge 0$, the subset ${\bf H}_m \subset {\bf H}$ satisfies ACC by
Lemma~\ref{ACC_for_AZ}.
 By Proposition~\ref{mod-gen},$(a)$, the set $\{ P(z)-h(z) | h(z) \in {\bf H}_m \}$  of
Hilbert series
 of quotient modules $P / L$ satisfies ACC, so, ${\bf H}_m$ satisfies DCC.
 By Koenig Lemma~\ref{Koenig}, this means that every set  ${\bf H}_m$ is finite.

 ${\sf (ES+EG)}  \Longrightarrow    {\sf (EC)}$.
For $d > 0$, let $P^d$ be the set of all submodules in $P$
generated in degrees at most $d$. Let $p(z)$ denote the
polynomial $\sum_{i \le d} z^d \dim P_i$. Every submodule
$L \in P^d$ is isomorphic to a quotient module of a free
module $M_L = X_L \otimes A$, where $X_L(z) \le p(z)$, by a
submodule $K_L$ generated by a minimal set of relations of
$L$. Since $P$ is coherent, $L$ is finitely presented, that
is, $K_L$ is finitely generated. The condition $X_L(z) \le
p(z)$ implies that there are only finite number of
possibilities for $X_L(z)$, so, there are only finite
number of isomorphism classes of $M_L$. Since $P$ is
effective for series, there are also only finite number of
possibilities for Hilbert series $L(z)$. Finite free module
$M_{L}$  is effective for generators, so, given a Hilbert
series  $K_L(z) = M_{L}(z) - L(z)$, there are only finite
number of possibilities for $m_1 (L) = m_0 (K_L)$. The set
of such Hilbert series $K_L(z)$ is finite, therefore, $P$
is effectively coherent.

${\sf (EC)}  \Longrightarrow    {\sf (ES)}$. Since $A$ is
effectively coherent, it is coherent, and we have $\dim
\tor_i^A(k,k) < \infty $ for all $i \ge 0$, hence we can
apply  to the set of Hilbert series of the quotient modules
$P/M$
 (which has the same cardinality as ${\bf H}_d$)
Proposition~\ref{mod-gen},~$(b)$
with $m= \dim H_0 P$, $p_2 = m(P)$, $q = d$, and
$r = D_P (d)$.
\end{proof}

\begin{cor}
\label{str_noet_effect}
Every strongly Noetherian algebra over an algebraically
closed field  is effectively Noetherian.
\end{cor}

We do not know, in what (geometrical?) terms the function $D(d)$
could be estimated for general strongly Noetherian algebras.

%It has been  shown in~\cite{az} that, if $P$  is a finite module over
%a   strongly Noetherian algebra over  an algebraically
%closed field, then the set of isomorphism classes of its submodules
%with given Hilbert series $h(z)$ is parameterized by a commutative projective %scheme.
%Therefore, we have

% \begin{cor}
%Let $k$ is algebraically closed, and let  $P$ be  a finite module over
%a   strongly Noetherian algebra $R$.
%Then the set of its submodules $L$ generated in degrees at most $d$
%is parameterized by a commutative  projective scheme.
%\end{cor}

\begin{quest}
\label{Q-PI-SN}
For a Noetherian PI algebra $A$ over an algebraically
closed field, are there estimates for the syzygy degree function
$D_A(d)$ in terms of its generators, relations,
and identities?
\end{quest}

A partial converse to the implication $(EC) \Longrightarrow (ES)$ of
Theorem~\ref{main_effect}
is given by

\begin{prop}
Let $A$ be a coherent algebra of global dimension 2.
Then  $A$ is effective for series if and only if it is effectively coherent.
\end{prop}

\begin{proof}

The ``if'' part is proved in Theorem~\ref{main_effect}; let
us proof the ``only if'' part.

Let $I$ be a right ideal of $A$ with $m(I) \le d$. Since
$I$ has projective dimension at most one, its minimal free
resolution has the form
$$
  0 \to V_2 \otimes A \to   V_1 \otimes A \to I \to 0,
$$
where $V_1, V_2$ are finite--dimensional vector spaces.
Given $d$, there are only finite number of possibilities
for the Hilbert series $V_1(z)$ of the space of generators
of $I$.  Since  $A$ is effective for series, there is also
only a finite number of possibilities for
$$
    V_2(z) = A(z)^{-1} \left( V_1(z) A(z) - I(z) \right) .
$$
Thus, there exists a constant $D = D(d)$ such that $m_1(I) = m(V_2) \le D$.
\end{proof}

However, in general an effective for series algebra of
global dimension two must not be effectively coherent,
e.~g., every finitely generated algebra over a finite field
is effective for series, but is not necessarily coherent
(for example, the algebra $k \langle x,y,z,t |  zy-tz, zx
\rangle $ has global dimension two but is not
coherent~\cite[Proposition~10]{pi2}).

Examples of Noetherian but not strongly Noetherian algebras
has been found by Rogalski in~\cite{rog}. The simplest
example is the ring $R_{p,q}$ generated by two Eulerian
derivatives with two parameters $p,q \in k^*$~\cite{jor}.

\begin{prop}
Assume that either $\char k = 0$ or $\td k \ge 2$. Then for
some $p,q \in k^*$ the  algebra $R_{p,q}$ is Noetherian but
neither effective for generators nor effective for series;
in particular,  $R_{p,q}$ is not effectively coherent.
\end{prop}

\begin{proof}

By~\cite[Example~12.8]{rog}, there are $p,q$ such that the
algebra  $R = R_{p,q}$ is Noetherian. To show that
$R_{p,q}$ is not effective in any sense, we will use some
properties of its point modules~\cite[Section~7]{jor}.
Recall that a {\it point module} is a cyclic module $M =
R/I$ with Hilbert series $M(z) = (1-z)^{- 1}$; the ideal
$I$ is called  a {\it point ideal}. It is shown
in~\cite[Section~7]{jor} that for every $d>0$ there are two
non-isomorphic point modules such that their truncations
$M/M_{<d}$ are isomorphic. In particular, a point ideal $I$
can have generators of arbitrary high degree. Since the
Hilbert series  $I(z) = h(z):= R(z) - (1-z)^{-1}$ is the
same for all point ideals, it follows that $R$ is not
effective for generators.

More precisely, there is a family of point modules $M(n,c )
= R/I(n,c ) $ (where $n\ge 0$ is an integer, $c \in {\PP}^1
(k)$), such that for $n \ge 2$  they have several relations
of degree at most 3 and one relation of degree
$(n+1)$~\cite[Section~7]{jor}. Let $I'(n,c) = I(n,c )_{\le
3}R$. Then $I'(n,c)(z) = h(z) + o(z^{n})$  while
$I'(n,c)(z) \ne h(z) + o(z^{n+1})$, therefore, the set of
Hilbert series of the ideals $I'(n,c)$  is infinite. Since
all these ideals are generated in degrees at most 3, $R$ is
not effective for series, whence is not effectively
coherent by Theorem~\ref{main_effect}.
\end{proof}

\section{Coherent families and universally coherent algebras}

\label{rate}

\subsection{Algebras with coherent families}

\begin{definitia}
\label{rate_fil_def}
\label{rate_def}
Let $R$ be a
%standard (i.e., degree-one generated)
finitely generated graded algebra, and let ${\bf F}$ be a
set of finitely generated right ideals in $R$. The family
${\bf F}$ is said to be coherent if:

1) zero ideal and the maximal homogeneous ideal $\overline
R$ belong to ${\bf F}$, and

2) for every $0 \ne I \in {\bf F}$
    there are ideal $I \ne J \in {\bf F}$ and a homogeneous element
    $x \in I$ such that $I = J + x R, m(J) \le m(I)$, and the ideal
    $N = (x:J) :=  \{ a \in R | xa \in J \}$ belongs to ${\bf F}$.
\end{definitia}

A coherent family ${\bf F}$ is said to be  of degree $d$,
if $m(I) \le d$ for all $I \in {\bf F}$,
i.e., all
its members are generated in degrees at most $d$.
%Koszul filtrations are exactly rate filtrations of degree one.

If $R$ is commutative and standard (i.e., degree-one
generated), coherent families are called   {\it generalized
Koszul filtrations}~\cite{cnr}. In particular, it is shown
in~\cite[Theorem~2.1]{cnr} that coordinate rings of certain
sets of points of  the projective space ${\bf P}^n$ admit
generalized Koszul filtrations of finite degrees. In the
non-commutative settings I cannot imagine a coherent family
as a filtration, that this new notion is introduced. The
term ``coherent family'' itself was proposed by L.
Positselski.

A coherent family of degree one is called {\it Koszul
filtration}: in this case, every ideal $I \in {\bf F}$ is
generated by linear forms. This concept has been introduced
for commutative algebras and investigated in several
papers~\cite{crv,ctv,i-kos,conca,con2}. In particular,
every coordinate ring of an algebraic curve   admits a
Koszul filtration provided that it is quadratic. The
non-commutative version of Koszul filtrations has been
considered in~\cite{kfx}; for example, every generic
$n$--generated quadratic algebra $R$ admits a Koszul
filtration if either it has less than $n$ relations or is
such that $\dim R_2 < n$.

Recall that a ring $R$ is called {\it (right) coherent} if
every map $M \to N$ of two finitely generated (right) free
$R$-modules has finitely generated kernel.  If the algebra
$R$ is graded, two versions of coherence can be considered,
``affine'' (general)  and ``projective'' (where all maps
and modules are assumed
 to be graded): the author does not know whether these concepts are
equivalent or not for connected algebras.
%The (projective) coherent rings may be considered as a basic of non-commutative
%geometry  instead of Noetherian rings~\cite{po}.

One of equivalent definitions of a coherent ring is as
follows~\cite[Theorem~2.2]{chase}: $R$ is coherent iff, for
every finitely generated ideal $J = JR$ and every element
$x \in R$, the ideal $N = (x:J)$  is finitely generated.
The similar criterion holds for projective coherence. This
definition is similar to our definition of coherent family,
as shows the following

\begin{prop}
\label{rate-coh}

For a standard algebra $R$, the following two statements
are equivalent:

(i) $R$ is projective coherent;

(ii) all finitely generated homogeneous ideals in $R$ form
a coherent family.
\end{prop}

\begin{proof}

The implication $(i) \Longrightarrow (ii) $ follows from
the  criterion above. The dual implication $(ii)
\Longrightarrow (i)$ follows from
Proposition~\ref{fil-rate},~$a)$.
\end{proof}

Note that the existence of a coherent family (even a Koszul
filtration) is not sufficient for an algebra to be
(projective) coherent. Indeed, the algebra  $A = k \langle
x,y,z,t |  zy-tz, zx \rangle $ admits a Koszul filtration
(it is initially Koszul with the Groebner flag $(x,y,t,z)$,
see~\cite[Section~5]{kfx}), but it is not projective
coherent, since the annihilator $\Ann_A z$ is not finitely
generated~\cite[Proposition 10]{pi2}.

The following property of coherent families gives linear
bounds for degrees of solutions of some linear equations.

Let us recall some notations.
% of~\cite{cnr}.
By definition, the {\it rate}~\cite{brate} of a (degree one generated) algebra $R$
is the number
$$
    \rate R = \sup_{i \ge 2} \{ \frac{m_i(R) - 1}{i-1} \}.
$$
For commutative standard algebras~\cite{an2,brate} as well as for
non-commutative algebras with finite Groebner basis of relations~\cite{an2}
 the rate is always finite.
The rate is equal to 1 if and only if $R$ is Koszul. If an
algebra has finite rate, then its Veronese subring of
sufficiently high order is Koszul~\cite{brate}.

The following Proposition (part $b)$) was originally proved
for commutative algebras in~\cite[Proposition~1.2]{cnr}. In
fact, it holds for non-commutative ones as well.

\begin{prop}
\label{fil-rate}

a) Let ${\bf F}$ be a coherent family in an algebra $R$.
Then the trivial $R$--module $k_A$ and every ideal $I \in {\bf F}$
have free resolutions of finite type.

b) Assume in addition, that
the  coherent family  ${\bf F}$ has degree at most $d$.
 Then
$$
      m_i(I)  \le m(I) + d i
$$
for all $i \ge 1$ and $I \in {\bf F}$.
In particular, if  $R$ is generated in degree one, then
$$
      \rate R \le d.
$$
\end{prop}

\begin{proof}

Following the arguments of~\cite[Proposition~1.2]{cnr}, we
proceed by induction in $i$ and in $I$ (by inclusion).
First, note that the degree $c$ of $x$ in
Definition~\ref{rate_fil_def} cannot be greater than
$m(I)$, and so, in the case~$b)$ it does not exceed $d$.
Taking $J$ and $N$ as in Definition~\ref{rate_fil_def}, we
get the exact sequence

$$
      0 \to J \to I \to R/N [-c] \to 0,
$$
which gives to the following fragments of the exact
sequence of $\tor$'s:
$$
    H_1 (J)_j \to H_1 (I)_j \to k_{j-c}
$$
and
$$
  H_i (J)_j \to H_{i} (I)_j \to H_{i-1} (N)_{j-c}
$$
for $i \ge 2$.

By induction, the first and the last terms in these triples
are finite-dimensional $k$--modules, so the middle one is.
This proves $a)$. Also, in the case $b)$, the first term
vanishes for $j \ge m(J)+di$, and the third term vanishes
for $j-c \ge m(N) + d(i-1)$. Since $m(J)\le m(I)$ and $m(N)
\le d$, they both vanish for all $j \ge m(I) + di$, so that
the middle term vanishes too.
\end{proof}

\subsection{Hilbert series and coherent families}

\begin{prop}
Let
${\bf F}$ be a coherent family of degree $d$ in an algebra $R$.
Then the set of all Hilbert series of ideals  $I \in {\bf F}$
is finite.
\end{prop}

\begin{proof}

By Proposition~\ref{fil-rate}, for every ideal $I \in {\bf
F}$ we have $m(I) \le d$  and $m_1(I) \le 2d$. Since
$\overline R \in {\bf F}$, we have also $m_1(R) \le d,
m_2(R) \le 2d$, and $m_3(R) \le 3d$. Thus, we can apply
Corollary~\ref{mod}.
\end{proof}

The following property of algebras with coherent families
seems to be the most interesting. Its analogue for Koszul
filtrations  has been proved in~\cite[Theorem~3.3]{kfx}.

\begin{theorema}
\label{rate ratio}
Suppose that an algebra $R$ has a
coherent family ${\bf F}$ of degree $d$.
Then $R$ has rational Hilbert series, as well as
every ideal $I \in {\bf F}$.

If the set ${\mathcal H}ilb$ of all Hilbert series of ideals
$I \in {\bf F}$
 contains at most $s$ nonzero elements, then the degrees
of numerators and denominators of these rational functions are not
greater than $ds$.
\end{theorema}

\begin{proof}

Let ${\mathcal H}ilb = \{I_0(z) = 0, I_1(z), \dots,
I_s(z)\}$, where $I_1, \dots, I_s$ are some nonzero ideals
in ${\bf F}$. By definition, for every  nonzero ideal $I =
I_i$ there are ideals $ J = J(I), N = N(I) \in {\bf F}$
such that $J \subset I, J \ne I$,     and for some positive
 $c = c(I)$ the following triple is exact:
$$
     0 \to J \to I \to R/N [-c] \to 0.
$$
Taking the Euler characteristics, we deduce
\begin{equation}
\label{rate_hilb}
     I(z) = J(I)(z) + z^c (R(z) - N(I)(z)).
\end{equation}
Let us put $J^{(1)} := J(I)$ and $J^{(n+1)} := J\left(
J^{(n)} \right)$. Since all ideals $J^{(n)}$  are generated
by subspaces of a finite-dimensional space $R_{\le d}$, the
chain
$$
    I \supset J^{(1)} \supset J^{(2)} \supset \dots
$$
contains only a finite number of nonzero terms. Applying
the formula~(\ref{rate_hilb}) to the ideals $J^{(n)}$, we
obtain a finite sum presentation
$$
    I(z) = z^{c(I)} (R(z) - N(I)(z))
            + z^{c(J^{(1)})} ((R(z) - N(J^{(1)})(z)))
            + z^{c(J^{(2)})} ((R(z) - N(J^{(2)})(z)))
        + \dots
$$
Thus
$$
    I_i(z) = \sum_{j=1}^s a_{ij} \left( R(z) - I_j(z) \right) + a_{i0}R(z),
$$
where $a_{ij} \in z {\bf Z} [z]$ (the last term corresponds
to the cases $N(J^{(t)}) = 0$).

Let $H = H(z)$ be the  column vector $[I_1(z), \dots,
I_s(z)]^{t}$, let $A$ be the matrix $(a_{ij}) \in M_s(z
{\bf Z}[z])$, let $H_0 = H_0(z)$ be the  column vector $[
a_{1,0}, \dots, a_{s,0}]^{t}$, and let $e$ be the unit
$s$--dimensional column vector. Then we have
$$
    H = A \left( R(z)e - H \right) + R(z) H_0,
$$
or
$$
   \left( A + E \right) H = R(z) (Ae +H_0),
$$
where $E$ is the unit matrix.

The determinant $D(z) = \det \left( A + E \right) \in {\bf
Z}[z]$ is a polynomial of degree at most $sd$. It is
invertible in ${\bf Q} \left[ [z] \right]$. Then $\left( A
+ E \right)^{-1} = D(z)^{-1} B$ with $B \in M_s({\bf
Z}[z])$. The elements of $B$ are $(s-1)\times (s-1)$ minors
of $\left( A + E \right)$, so, their degrees do not exceed
$d(s-1)$. We have
$$
    H = R(z) D(z)^{-1} C,
$$
where $C = B(Ae +H_0)  \in z {\bf Z}[z]^s$. So, for every
$1 \le i \le s$ we have
$$
    I_i(z) = R(z) C_i(z) D(z)^{-1}.
$$

Assume that $\overline R = I_s$. Then
$$
   \overline R (z) = R(z) -1 = R(z) C_s(z) D(z)^{-1},
$$
therefore,
$$
     R(z) = \frac{D(z)}{ D(z) - C_s(z)},
$$
so, $R(z)$ is a quotient of two polynomials of degrees at most $sd$.
By the above,
$$
    I_i(z) = R(z) C_i(z) D(z)^{-1} = \frac{ C_i(z)}{ D(z) - C_s(z)},
$$
so, the same is true for the Hilbert series $I_i(z)$.
\end{proof}

\subsection{Universally coherent algebras}

An algebra $R$ is called {\it universally $d$--coherent} if
all ideals in $R$ generated in degrees at most $d$ form a
coherent family.  Universally 1--coherent algebras are
called {\it universally Koszul} (because all their ideals
generated by linear forms are Koszul modules). Commutative
universally Koszul  algebras have been considered
in~\cite{conca, con2}. In particular, commutative monomial
universally Koszul algebras are completely classified
in~\cite{con2}.

The following criterion is a consequence of
Proposition~\ref{fil-rate}.

\begin{prop}
\label{D<t+d}
The following two conditions for a connected algebra $R$ are equivalent:

(i) $R$ is  universally $d$--coherent;

(ii) $D_R(t) \le t+d$ for all $t \le d$.
\end{prop}

\begin{proof}

$(i) \Longrightarrow (ii)$.

By Proposition~\ref{fil-rate}, b), we have $m_1(I) \le m(I)
+ d$ if $t = m(I) \le d$.

$(ii) \Longrightarrow (i)$.

By $(ii)$, we have  $m_1(I) \le m(I) + d$ if $t = m(I) \le
d$. Let $x = \{ x_1, \dots, x_r \}$ be a minimal system of
generators for $I$ with $\deg x_1 = t$, let $J = x_2 R +
\dots + x_r R$, and let $N = N(I) = (x:J)$. Obviously, $N$
is the (shifted by $t$) projection of the syzygy module
$\Omega = \Omega (x_1, \dots, x_r)$ onto the first
component, so that $m(N) + t \le m(\Omega) = m_1(I) \le
t+d$. Thus, $m(N (I)) \le d$ for every $I$ with $ m(I) \le
d $. By Definition~\ref{rate_fil_def}, this means that all
such ideals $I$ form a coherent family.
\end{proof}

\begin{corollarium}
\label{D<2d}
The following two conditions for a connected algebra $R$
are equivalent:

(i) $R$ is universally $d$--coherent for all $d \gg 0$;

(ii)$D_R(d) \le 2d$ for all $d \gg 0$.

In particular,  in this case $R$ is effectively coherent.
\end{corollarium}

We call an algebra $R$ satisfying either of the equivalent
conditions of this Corollary~\ref{D<2d} simply {\it
universally coherent}.

Every   generated in degree one universally $d$--coherent
algebra is a finite Koszul module over a universally Koszul
algebra, namely, over its Veronese subalgebra. This follows
from

\begin{prop}
Let $R$  be a universally $d$--coherent algebra generated in degree one.
Then its Veronese subalgebra $R^{(d)}$ (with grading divided by $d$)
is universally Koszul.
\end{prop}

\begin{proof}

Let $B$ be the algebra  $R^{(d)}$ with grading divided by $d$.
By Proposition~\ref{D<t+d}, we have to check that $D_B (1) \le 2$.
By Propositions~\ref{DVerona} and~\ref{cap-ann}, we have
$$
     D_B^\cap (1) \le d^{-1} D_{R^{(d)}}^\cap (d) \le
     d^{-1} ( D_R ^\cap (d) + d-1) \le  3 -d^{-1},
$$
and, analogously,
$$
     D_B^{\ann} (t) \le   3 -d^{-1}.
$$
By Proposition~\ref{cap-ann}, we have
$$
      D_B (1) \le  \max \{ D_B^\cap (1), D_B^{\ann} (1)  \} \le 2.
$$
\end{proof}

In particular, in a universally coherent algebra all
Veronese subrings of sufficiently high order are
universally Koszul, thus the algebra itself is a finite
direct sum of Koszul modules over every such Veronese
subalgebra.

\begin{prop}
Every finitely presented module over a  universally  coherent algebra has rational
Hilbert series.
\end{prop}

\begin{proof}
Let $R$ be a universally coherent algebra. Since every
finitely generated ideal in $R$ is a member of a coherent
family, it has rational Hilbert series, as well as the
algebra $R$ itself.

Let $M$ be a finitely presented $R$--module, and let
$$
    0 \to K \to F \to M \to 0
$$
be its finite presentation with free $F$. Since $M(z) = F(z) - K(z)$, it is sufficient to
show that
the finitely generated submodule $K$ of a free module $F$ has rational Hilbert series.

Let $d = m(K)$, and let $K = L + xR$,  where $\deg x = d$.
By the induction in the number of generators of $K$, we can
assume that $L(z)$ is rational.

Let $J = (x : L)$. Since $R$ is coherent,
this ideal is finitely generated, therefore, its Hilbert series is rational.
The exact triple
$$
      0 \to L \to K \to R/J [-d] \to 0
$$
gives the formula $K(z) = L(z) + z^d (R(z) - J(z))$.
Thus $K(z)$ is rational.
\end{proof}

\subsection{Examples: monomial algebras and their generalizations}

\label{perera}

Universal coherence is the strongest property among all
considered above. It is a fortunate surprise that some
interesting classes of noncommutative algebras are indeed
universally coherent.

All finitely presented algebras whose relations are {\it
monomials} on generators are universally coherent (see
below). It is an important point in non-commutative
computer algebra that many properties of finitely presented
monomial algebras are inherited by the algebras with finite
Groebner bases (for example, they have finite rate, and
quadratic ones are necessary Koszul). However, the latter
are not in general coherent, see examples in~\cite{pi2};
the reason is that while the algebra itself has finite
Groebner basis, a finitely generated one-sided ideal in it
sometimes can have only an infinite one (in contrast to the
monomial case~\cite[Theorem~1]{pi1}). That is why a new
class of algebras between these two has been
introduced~\cite{pi2}, so-called {\it algebras with
$r$--processing}.

In general, these algebras are not assumed to be graded but
still finitely generated. Let $A$ be a quotient algebra of
a free algebra $F$ by a two-sided ideal $I$ with a Groebner
basis $G = \{ g_1, \dots, g_s \}$. For every element $f \in
F$, there is a well-defined normal form $N(f)$  of $f$ with
respect to $G$~\cite[Subsection~2.3]{ufn} .

The algebra $A$ is called {\it algebra of $r$--processing}
for some $r \ge 0$, if for any pair $p,q \in F$ of normal
monomials, where $q = q_1 q_2, \deg q_1 \le r$, we have
$$
      N(p q) = N(p q_1) q_2.
$$
The simplest example is an algebra $A$ whose relations are
monomials of degree at most $r+1$.

A simple sufficient condition for an algebra to have this
property is as follows. Consider a graph $\Gamma$ with $s$
vertices marked by $g_1, \dots, g_s$ such that an arrow
$g_i \to g_j$ exists iff there is an overlap between any
{\it non-leading} term of $g_i$  and {\it leading} term of
$g_j$. If $\Gamma$ is acyclic, then  $A$ is an algebra with
$r$--processing for some $r$. The simplest case is when the
monomials in the decompositions of the relations of the
algebra do not overlap each other. See~\cite{pi2} for other
sufficient conditions and for a way how to calculate the
number $r$ for given $\Gamma$.

The main property of algebras with $r$--processing is that every
finitely generated right ideal in such algebra has finite Groebner basis.
In the case of standard degree--lexicographical order on monomials,
the degrees of its elements do not exceed $m(I) + r -1$.
Moreover, the degrees of relations of $I$ do not  exceed the number $m(I) + 2r-1$; in particular, they are coherent.
Therefore, we have

\begin{prop}
Let $R$ be a connected algebra with $r$--processing.
Then $D_R(d) \le d + 2r$. In particular, it is universally coherent.
\end{prop}

For example, any  algebra with $1$--processing is
universally Koszul (algebras with $1$--processing were
separately considered in~\cite{iou}). So, all quadratic
monomial algebras are universally Koszul (unlike the
commutative case, see~\cite{con2}).

\end{document}